\documentclass{ifacconf}

\usepackage{natbib} 
\usepackage{graphicx}

\usepackage{enumitem}
\usepackage{amsmath}
\usepackage{amssymb}
\usepackage{amsfonts}
\usepackage[dvipsnames]{xcolor}

\usepackage{mathtools}
\mathtoolsset{showonlyrefs}

\newtheorem{theorem}{Theorem}[section]
\newtheorem{corollary}[theorem]{Corollary}
\newtheorem{proposition}[theorem]{Proposition}
\newtheorem{lemma}[theorem]{Lemma}

\theoremstyle{definition}
\newtheorem{remark}[theorem]{Remark}
\newtheorem{assumption}[theorem]{Assumption}

\newcommand{\R}{\mathbb{R}}
\newcommand{\N}{\mathbb{N}}
\newcommand{\diff}{\, \mathrm{d}}
\newcommand{\id}{\mathrm{id}}

\newcommand{\E}{\mathbb{E}}

\newcommand{\Yavg}{Y_{\mathrm{avg}}}
\newcommand{\Yb}{\mathbf{Y}}

\newcommand{\cst}{\mathrm{Cst}}
\newcommand{\muprior}{\mu_{\mathrm{prior}}}
\newcommand{\prior}{\mathrm{prior}}
\newcommand{\post}{\mathrm{post}}
\newcommand{\Tr}{\mathrm{Tr}}
\newcommand{\U}{\mathcal{U}}

\DeclareMathOperator{\rank}{rank}
\newcommand{\J}{\mathcal{J}}

\newcommand{\uopt}{\hat u^*}

\DeclareMathOperator*{\argmin}{arg\,min}

\begin{document}
\begin{frontmatter}

\title{A case study in ensemble optimal control for Bayesian input design}
\author[Ludo]{L. Sacchelli} 
\author[Ale]{A. Scagliotti} 

\address[Ludo]{Inria, Université Côte d’Azur, CNRS, LJAD, France}
\address[Ale]{
Department of Mathematics, CIT School, Technical University of Munich, Germany, and Munich Center for Machine Learning}

\thanks{emails: \texttt{ludovic.sacchelli@inria.fr, scag@ma.tum.de}}


\begin{abstract}
We discuss the problem of input design for uncertainty reduction in a parameter estimation procedure. Assuming a linear continuous-time control system with noisy measurements, we formulate an objective of variance reduction in a Bayesian Gaussian setting as an optimal control problem and analyze it from a geometric control perspective.  The resulting cost functional depends on the unknown parameter, we compare the optimal control approach with a non-standard alternative inspired by ensemble control, where the cost is averaged over the prior distribution after computation, rather than before. This requires the statement of a generalized Pontryagin's maximum principle adapted to Gaussian distributions.
\end{abstract}
\end{frontmatter}

\section{Introduction}
In the context of system identification, experiment design aims to plan data acquisition so as to make parameter estimation as informative and efficient as possible. This is achieved by formulating and solving optimization problems that balance the reduction of parameter uncertainty with the cost or burden of experimentation. See, e.g., \cite{pukelsheim2006optimal,fedorov2013theory}.
In particular, for input–output systems, the problem often becomes one of input design: selecting an input that yields the most informative data for parameter identification. By coupling the experiment with the identification procedure, one can quantify the improvement achieved by each new experiment.
A common approach relies on the Fisher information matrix, which one seeks to maximize \cite{pronzato2008optimal}. An alternative, equally popular approach adopts a Bayesian perspective, aiming instead to minimize, according to a chosen metric, the variance of the posterior distribution \cite{ryan2016review}.

A central difficulty in experiment design for identification procedures lies in the fact that, when estimating a parameter, the outcome of the experiment generally depends on the very parameter being estimated. In other words, for a system depending on a parameter $\theta$, the optimization objective function $J$ associated with the input $u$ will itself depend, possibly implicitly, on $\theta$, i.e., $J(u,\theta)$. This issue was recognized early in the development of the theory. The standard approach is to replace the unknown value of $\theta$ with a nominal value $\bar\theta$, typically an a priori estimate.
It was soon understood (\cite{Hjalmarsson1996}), however, that iterative strategies offer a natural improvement: by repeating experiments, one can refine the parameter estimate successively, using the result of experiment number $n$ to design experiment number $n+1$.

Nevertheless, nominal parameter values fail to account for uncertainty in the prior estimate.
To achieve more robust designs, one can explicitly incorporate this uncertainty into the procedure
(\cite{bombois2021robust,petsagkourakis2021safe}). 
In this respect, Bayesian formulations are particularly well suited. When a prior estimate $\theta_\prior$, distributed according to a prior measure $\mu_\prior$, is available, the classical approach consists in optimizing the cost 
$J(u,\E_{\mu_\prior}[\theta])$,
where the unknown parameter is replaced by its prior mean.
A more robust alternative is to optimize the expected cost
$\E_{\mu_\prior}[J(u,\theta)]$,
thereby accounting for parameter uncertainty at the design stage.

From the systems viewpoint, the optimal control approach to experiment design has received some attention. Linear cases, considered in the present work, were explored early on (\cite{Kalaba1982}), often reducing to variants of LQR control. Nonlinear studies remain sparse and are mostly numerical, relying on either MPC-based (\cite{qian2017optimal}) or shooting-based methods (\cite{bock2013model}). In both these settings and ours, relying on a nominal value for the unknown parameter leads to a standard control problem. In contrast, optimizing the expected cost
is more involved: it requires considering the input-output system as a distributed system with respect to the parameter, and then to compute an averaged cost according to the prior distribution. This framework, however, falls within optimal ensemble control. Ensemble control is an emerging field of control theory that allows to recast uncertainty in parameters as control of infinite dimensional systems. While it has seen some initial interest in the communities of quantum and multi-agent control, recent work has been devoted to the extension of Pontryagin's maximum principle (PMP) to ensemble control addressing averaged cost (\cite{bettiol2019necessary,scagliotti2023optimal}).

In this paper, we examine optimal control formulations of an input design problem. We consider a linear system in a Gaussian Bayesian estimation setting and formulate an optimization problem aimed at reducing posterior uncertainty while controlling state magnitude. The problem is first analyzed in its classical form, then extended to the averaged ensemble control framework. Existing results on the ensemble PMP do not apply to linearly parameterized systems with Gaussian-distributed parameters. We therefore introduce technical arguments that generalize the ensemble PMP to this specific structure.

\section{Input design problem}

\subsection{Setting}

Let $n,m,p,q\in \N$. We set $A\in \R^{n\times n}$, $B\colon \R^p\to \R^{n\times m}$, $C\in \R^{q\times n}$, and $\sigma\in \R^{q\times q}$. We assume the mapping $B\colon \R^p\to\R^{n\times m} $ to be affine, of the form
$$
B^\theta=B_0+\theta_1 B_1 + \dots +\theta_p B_p, \quad B_0,\dots , B_p\in \R^{n\times m}.
$$
Let $T>0$ and for now let $u\in L^\infty([0,T],\R^m)$ (for the presumed Lebesgue measure).
For a given parameter $\theta\in \R^p$,  we consider the controlled dynamical system, initiated at $x(0)=0$:
\begin{equation}\label{E:sde1}
\begin{cases}
\dot x^\theta=A x^\theta + B^\theta u,
\\
y=C x^\theta.
\end{cases}
\end{equation}
The notation $x^\theta$ highlights the implicit dependence of the trajectory to $\theta$.
Consider $(B^\theta)^i\in \R^n$ the $i$-th column of $B^\theta$, $1\leq i \leq m$. Then we define $B^i$ the $n\times p$ matrix constituted of the $i$-the columns of $B_1,\dots , B_p$, (and $B_0^i$ the $i$-th column of $B_0$) so that $(B^\theta)^i=B_0^i+B^i\theta$.
\begin{assumption}\label{A:Btheta}
We assume that $q\leq p\leq n$, $m\leq n$ (less measurements than parameters). Furthermore, 
the family $\left(  B^i\right)_{1\leq i\leq m}$ is free (no control redundancy).
\end{assumption}

\begin{remark}
In the optimal control analysis we provide, we make the assumption $q=1$, the analysis prior to that point holds for arbitrary $q$ however.
\end{remark}

We assume that measurements are noisy and we design an estimator around an averaging of the noise. That is, we consider the equation (with $W(t)$ a standard $q$-dimensional Brownian motion, $\sigma\in \R^q\times \R^q$ of rank $q$)
\begin{equation} \label{eq:Y_signal}
\Yavg=\frac{1}{T}\int_0^T y(t)\diff t + \frac{\sigma}{T}\int_0^T  \diff W(t).
\end{equation}
We wish to use this expression to estimate $\theta$. We have
\begin{equation}
\Yavg
=
Y_0+\sum_{i=1}^p\theta_j Y_j+\varepsilon
\end{equation}
where the noise $\varepsilon$ follows a centered Gaussian law of variance $\frac{1}{T}\sigma\sigma^\top$ (i.e. $\varepsilon\sim \mathcal{N}\left(0,\frac{1}{T}\sigma\sigma^\top\right)$) and 
\begin{equation}
Y_j=\frac{1}{T}\int_0^T\int_0^tC\e^{A(t-s)}B_ju(s)\diff s \diff t,\quad j=0,\dots ,p.
\end{equation}
Notice that if we introduce the functions
\begin{equation}
\psi_j(s)=\int_s^T B_j^\top \e^{A^\top(t-s)}C^\top\diff t, \quad j=0,\dots ,p,
\end{equation}
then, after some Fubini manipulation, we get
\begin{equation}
Y_j=\frac{1}{T}\int_0^T \psi_j(t)^\top u(t)\diff t,\qquad j=0,\dots ,p.
\end{equation}
Finally, we introduce $\Yb\in \R^{q\times p}$ whose columns are the vectors $(Y_j)$, so that 
$
\Yb=\left(
\begin{array}{c|c|c}
Y_1&\cdots&Y_p
\end{array}
\right).
$
In particular
\begin{equation}
\Yavg
=
Y_0+\Yb\theta+\varepsilon.
\end{equation}

\subsection{The Bayesian input design problem}
\label{S:Bayesian}

Given a Gaussian prior distribution on the parameter $\theta\sim \mathcal{N}(\theta_\prior,\Sigma_\prior)$, we ask what would be the ideal choice of $u$ to minimize the variance of the Gaussian posterior.
For this we proceed to a Bayesian analysis of the experiment. Denoting by $\rho$ the associated likelihood functions, and with $\cst$ constants  not relevant to the analysis,
\begin{multline}
-2\log 
\rho
\left(
\Yavg\mid \theta
\right)
=\\
\left(\Yavg-Y_0-\Yb\theta\right)^\top T(\sigma\sigma^\top)^{-1}
\left(\Yavg-Y_0-\Yb\theta\right) 
+
\cst,
\end{multline}
while 
\begin{equation}
-2\log 
\rho
\left(
\theta
\right)
=
(\theta-\theta_\prior)\Sigma_\prior^{-1}(\theta-\theta_\prior)
+
\cst.
\end{equation}
As a result
\begin{equation}
\begin{aligned}
&-2\log 
\rho
\left(
\theta\mid\Yavg
\right)
=
(\theta-\theta_\prior)\Sigma_\prior^{-1}(\theta-\theta_\prior)+
\\&\hphantom{= }\;
\left(\Yb\theta+Y_0-\Yavg\right)^\top T(\sigma\sigma^\top)^{-1}
\left(\Yb\theta+Y_0-\Yavg\right)
+
\cst
\\&
=
\theta^\top
\left(
T\Yb^\top(\sigma\sigma^\top)^{-1}\Yb+\Sigma_\prior^{-1}
\right)\theta
\\&
\hphantom{= }\;
-
2\theta^\top\left(
T\Yb^\top (\sigma\sigma^\top)^{-1} (\Yavg-Y_0)
+
\Sigma_\prior^{-1}\theta_\prior
\right)
+\cst.
\end{aligned}
\end{equation}
Letting $S$ be the symmetric positive definite matrix such that $S^2=T (\sigma\sigma^\top)^{-1}$ ($\sigma$ has rank $q$),
we induce that the posterior must have
variance 
\begin{equation}
\Sigma_\post
=
\left(
\Yb^\top S^2\Yb+\Sigma_\prior^{-1}
\right)^{-1}
\end{equation}
and 
expected value
\begin{equation}
\theta_\post=
\Sigma_\post\left(
\Yb^\top S^2 (\Yavg-Y_0)
+
\Sigma_\prior^{-1}\theta_\prior
\right).
\end{equation} 
In the case where no prior knowledge is assumed, $\Sigma_\prior^{-1}$ and $\Sigma_\prior^{-1}\theta_\prior$ terms are omitted in the expressions above, but we do not discuss this case in the following.
Denoting by $\Sigma_\prior^{1/2}$ the square root of the positive definite  $\Sigma_\prior$, we have, as with any linear Gaussian Bayesian problem, 
\begin{equation}
\begin{aligned}
\Sigma_\post
&=
\Sigma_\prior
\left(
\id+
\Sigma_\prior^{1/2}\Yb^\top S^2\Yb\Sigma_\prior^{1/2}
\right)^{-1}
\preceq\Sigma_\prior.
\end{aligned}
\end{equation}
(Denoting by $\preceq$ Loewner's order.)

We now know enough to discuss input design objectives.
The goal is to lower uncertainty as much as possible so a standard criterion is D-optimality, where the objective is to design the control that minimizes $\log \det \Sigma_\post$.
This objective is relaxed using concavity of $\log \det$ into 
\begin{equation}\label{E:costTr}
\Tr
\left(
\Sigma_\prior^{1/2}\Yb^\top S^2 \Yb\Sigma_\prior^{1/2}
\right)
\longrightarrow \max.
\end{equation}
Recognizing the Frobenius inner product in \eqref{E:costTr}, this objective function is quadratic in $\Yb$, itself linear in $u$, thus being a quadratic objective function in $u$ (but not necessarily coercive). 
This relaxation is a frequent strategy to avoid the high non-linearity of D-optimality (although determinant is still polynomial) and recover quadratic objectives which have been explored via standard LQR techniques. 
In the present paper, we follow a slightly alternative approach to arrive at a similar simplification. Alternative to D-optimality, the E-optimality criterion aims at 
minimizing the largest eigenvalue of $\Sigma_\post$ (hence minimizing the spectral norm). This is equivalent to maximizing 
$
\sigma_{\min}\left(\Sigma_\prior^{-1}+\Yb^\top S^2 \Yb\right)
$.
The lowest eigenvalue, as a map over symmetric matrices, is not differentiable, so this objective presents its own difficulty.
However, for any $\vartheta\in \R^p$, we find $\vartheta^\top\Sigma_\post^{-1}\vartheta = \vartheta^\top\Sigma_\prior^{-1}\vartheta+(S\Yb \vartheta)^\top (S\Yb \vartheta)$. In particular, the rank of $\Yb^\top S^2\Yb$ is at most $q$. 

Assume $q=1$ from now on.  We formulate an objective that prioritizes increasing the smallest eigenvalue of $\Sigma_\prior^{-1}$. Let 
\begin{equation}
    V
    =
    \argmin_{\substack{X^\top X=1
    \\
    X\in \R^{p}}}
    X^\top \Sigma_\prior^{-1} X.
\end{equation}
\begin{remark}
The vector $V$ is not unique if the minimal eigenvalue has multiplicity. 
Multiple eigenvalues are non-generic for symmetric matrices (an open and dense subset of symmetric matrices have only simple eigenvalues), so we do not consider this case in the present analysis.
\end{remark}




We deduce the relaxed objective: align $S\Yb$ as much as possible with the subspace generated $V$, which is quantified via
$|S\Yb V|^2$.
This is again a quadratic objective. Since we assume $q=1$, $S\Yb V$ is a single real value which we rewrite as a scalar product. Here $S=\sqrt{T}\sigma$ and 
$\Yb=
\begin{pmatrix}
Y_1&\cdots&Y_p
\end{pmatrix}\in \R^{1\times p}$, where  $Y_j=\frac{1}{T}\int_0^T \psi_j(t)^\top u(t)\diff t$, $1\leq j\leq p$. We then have
\begin{equation}
\begin{aligned}
S\Yb V
&=
S
\sum_{j=1}^p
 Y_j V_j 
=
\frac{1}{T}
\sum_{j=1}^p
\int_0^T
 S V_j \psi_j(t)^\top u(s)
\diff s
\\
&=:
\frac{1}{T}
\int_0^T
\psi(t)^\top u(s)
\diff s,
\end{aligned}
\end{equation}
with $\psi=S \sum_{j=1}^p \psi_j V_j \in C^\infty([0,T],\R^m)$.
By sign symmetry, we can remove the square to arrive at the most reduced objective, that is 
\begin{equation}\label{E:costScal}
    S\Yb V
    \longrightarrow \max.
\end{equation}

\subsection{Optimal control point of view}
\label{S:optCon}

As an optimization problem, maximizing the cost function $S\Yb V=\langle\psi,u\rangle_{L^2([0,T],\R^m)}$
is ill-posed as a linear form in $u$ (the same goes if we keep the square). So we consider the problem under constraints on the $L^\infty$-norm of $u$. 
Assuming $\|u\|_{L^\infty}\leq m $, the optimal control is such that $u_i=m\,\mathrm{sign} (\psi^i)$. 
This control maximizes $\Yb$ along a specific direction within the given constraint, reflecting minimization of the influence of noise on the total measurement.  However these inputs do not regulate the state, which can become large, especially with large $T$ and $m$.
For this reason, we introduce a penalized objective function
\begin{equation}\label{eq:def_J_1}
J_\alpha(u,\theta)=\frac{1}{T}\int_0^T \psi(t)^\top u(t)\diff t
-
\frac{\alpha}{T}\int_0^T |x^\theta(t)|^{2}\diff t.
\end{equation} 
This cost depends on three parameters: $T,m,\alpha$. Because the initial condition is $x(0)=0$, $x$ is linear in u (that we denote $x_u$). Then for any $\lambda>0$, 
\begin{equation}
J_\alpha(\lambda u,  \theta)    
=
\langle\psi,\lambda u\rangle_{L^2}
+
\alpha \langle x_{\lambda u},x_{\lambda u}\rangle_{L^2}
=
\lambda J_{\lambda\alpha} (u,\theta).
\end{equation}
As a consequence, optimizing $J_\alpha$ under  $\|u\|_{L^\infty}\leq m$ is equivalent to optimizing $J_{\lambda\alpha}$ under $\|u\|_{L^\infty}\leq \lambda m$. For this reason we eliminate the role of $m$ by setting it equal to $1$. From now on, we also simply write $J$ instead of $J_\alpha$.

This optimization problem illustrates a fundamental challenge in input design: the cost functional $J$ depends on the unknown parameter $\theta$. A standard approach is to replace $\theta$ with a fixed estimate, typically obtained from prior knowledge or a preliminary identification step. In our Bayesian framework, this estimate is the prior mean $\theta_\prior$, leading to the following optimal control problem:
\begin{equation} \label{eq:problemClassic}
\begin{cases}
\dot x=Ax + B^{\theta_\prior} u, \quad x(0)=0,
\\
\|u\|_{L^\infty(0,T)}\leq 1,
\\
J(u,\theta_\prior)\longrightarrow \max.
\end{cases}
\end{equation}

In the present paper, we wish to discuss 
this classical formulation in comparison with a less classical one. If the prior estimate $\theta_\prior$ is inaccurate, the resulting input may be far from optimal. To mitigate this, we account for the full prior distribution rather than a single point estimate. Specifically, we reformulate the problem as an ensemble control problem, where each $\theta\in \R$ defines a system $x^\theta$
governed by the same input $u$. Instead of maximizing the cost for the expected parameter, we maximize the expected cost under the prior distribution $\muprior$, yielding the ensemble optimal control problem:
\begin{equation} \label{eq:problemEnsemble}
\begin{cases}
\dot x^\theta =Ax^\theta  + B^\theta u, \quad x^\theta (0)=0,
\\
\|u\|_{L^\infty(0,T)}\leq 1,
\\
\E_{\muprior}\left[J(u,\theta)\right]\longrightarrow \max.
\end{cases}
\end{equation}

The linearity of the problem imply symmetries that allow the ensemble formulation to be reduced to a (albeit higher dimensional) classical one. These symmetries do not extend to nonlinear cases which we intend to further explore, still,  they provide a useful ground truth for comparison with numerical ensemble methods.

\begin{remark} \label{rmk:exact_ensemble_comput}
For $0\leq i\leq p$, for $u\in L^\infty(0,T)$, let $z_i(u,\cdot):[0,T]\to \R^n$ be the absolutely continuous function defined by $z_i(u,t)=\int_0^t \e^{A(t-s)}B_i u(s)\diff s$.  Letting $\theta_0=1$ and $\tilde\theta=(1,\theta)\in \R^{p+1}$,  $x^\theta(t)=\sum_{i=0}^p \theta_i z_i(u,t)$. We can rewrite this as
\begin{equation}
x ^\theta(t) = 
\begin{pmatrix}
    \theta_0 \id 
    &
    \cdots
    &
    \theta_p \id
\end{pmatrix}
z(u,t)
\end{equation}
by stacking $(z_i)_{0\leq i\leq p}$ as a column. 
\begin{equation}
Q(\theta)=\begin{pmatrix}
    \theta_0 \id_{\R^n}
    \\
    \vdots
    \\
    \theta_p \id_{\R^n}
\end{pmatrix}
\begin{pmatrix}
    \theta_0 \id_{\R^n} 
    &
    \cdots
    &
    \theta_p \id_{\R^n}
\end{pmatrix}
=\tilde \theta \tilde \theta^\top
\otimes
\id_{\R^n},
\end{equation}
so that
$
x ^\theta(t)^\top x ^\theta(t)
=
z(u,t)^\top
Q(\theta)
z(u,t)
$
and
\begin{equation}
\E_{\mu_\prior}
\left[
|x^\theta(t)|^{2}
\right]
=
z(u,t)^\top
\E_{\mu_\prior}[Q(\theta)]
z(u,t).
\end{equation}
In particular, $\E_{\mu_\prior}[Q(\theta)]$ is a $n(p+1)\times n(p+1)$ symmetric block matrix, with top left block $\id_{\R^n}$, and bottom right block $\left(\Sigma_\prior+\theta_\prior \theta_\prior^\top\right)
\otimes
\id_{\R^n} \in \R^{np\times np}$, $0$ elsewhere.
The ensemble problem is thus equivalent to a deterministic $n(p+1)$-dimensional problem with $\tilde{A}=\mathrm{diag}(A,\dots, A)$, $\tilde{B}^\top =
\begin{array}{c|c|c}
(B_0^\top&\cdots&B_p^\top)
\end{array}$, $\dot z = \tilde A z + \tilde B u$. Then the quantity $\langle x^\theta,x^\theta\rangle_{L^2([0,T],\R^n)}=\langle z,Q(\theta)z\rangle_{L^2([0,T],\R^{n(p+1)})}$, and computing an averaged cost is just considering $\E_{\mu_\prior}[Q(\theta)]$.
\end{remark}

\section{Classical optimal control}\label{S:classic}

We use optimal control theory to derive optimality conditions characterizing the structure of optimal trajectories. This analysis is conducted from two viewpoints introduced in Section~\ref{S:optCon}: the classical formulation \eqref{eq:problemClassic}, addressed here, and the ensemble control strategy \eqref{eq:problemEnsemble}, treated next. Both flows are expected to share similar features: the $L^\infty$ constraint induces a bang-bang structure, possibly linked by singular arcs. Second-order conditions confirm their occurrence in the classical case. For the ensemble case, a dedicated version of Pontryagin’s maximum principle is required, and will be introduced in the next section.



Before their shape, let us justify the existence of solutions to \eqref{eq:problemClassic}. Similar ideas will later allow to prove existence of the solution of \eqref{eq:problemEnsemble} (see Lemma~\ref{T:existence_generalized}).
We begin with a convergence result for the trajectories of \eqref{E:sde1} generated by weakly-convergent controls. From general results that hold for control affine systems (see, e.g., \cite[Lemma~2.4]{scagliotti2023optimal}), we have the following.
\begin{lemma} \label{lem:conv_traj}
Let $(u_k)_k$ be a sequence such that $u_k 
\to u$ in the $L^\infty$ weak-$*$ topology\footnote{Recall that $u_k 
\to u$ in the $L^\infty$ weak-$*$ topology if we have convergence against $L^1$ test functions.}.
Then 
$$
\lim_{k\to \infty} \| x_{u_k}^\theta - x_u^\theta \|_{C^0([0,T],\R^n)} =0.
$$
\end{lemma}

\begin{proposition}
Problem \eqref{eq:problemClassic} admits a solution.
\end{proposition}
\begin{pf}
Let us consider a maximizing sequence $(u_k)_k$ \eqref{eq:problemClassic} and, recalling that $\|u_k\|_{L^\infty}<1$, we consider a (not relabeled) subsequence such that $u_k \to  u_\infty$ ${L^\infty}$  $*$-weakly.
Recalling that $J(u,\theta_\prior) =  \langle \psi_1, u\rangle_{L^2}  - \alpha  \langle x_u^{\theta_\prior}, x_u^{\theta_\prior} \rangle_{L^2}$, by virtue of Lemma~\ref{lem:conv_traj}, we obtain that 
\begin{equation*}
    \max_{\|u\|_{L^\infty}\leq 1} J(u,\theta_\prior) = \lim_{k\to \infty} J(u_k,\theta_\prior) =
    J(u_\infty,\theta_\prior).
\end{equation*}
\end{pf}

\subsection{Description of the optimal flow}\label{S:classical}

For this we assume that $\theta=\bar\theta\in \R^p$ is not distributed but a chosen beforehand single value, typically $\theta=\theta_\prior=\E_{\muprior}[\theta]$. We describe the solution of \eqref{eq:problemClassic} using Pontryagin's maximum principle (see, e.g., \cite{Trelat2024}).
We simplify the notations denoting $\bar B=B^{\bar \theta }$, $\bar x=x^{\bar \theta}$. 

\begin{assumption}
We assume $\rank(\bar B)=m$.
\end{assumption}
\begin{remark}
The above assumption does not necessarily hold for all $\bar\theta$ (e.g., when $\bar\theta=0$ and $\rank(B)<m$). However, if $\rank(B^\theta)=m$ for at least one $\theta$, it is then true for an open and dense subset of $\R^p$. 
\end{remark}


Let 
$H(t,\bar x, \bar p,p_0,u)=\bar pA\bar x + p_0 \alpha |\bar x|^2 + (\bar p \bar B-p_0\psi)u.$
($p_0\leq 0$ constant such that $(\bar p,p_0)$ is nontrivial).
Any optimal trajectory is the projection onto $\R^n$ of a curve $t\mapsto (\bar x(t),\bar p(t))$ in  $T^*\R^n$, solution to the Hamiltonian system associated with $H$ and with a control $u$ that maximizes $u\mapsto H(t,\bar x, \bar p,p_0,u)$ at all time.

First, the Hamiltonian system is 
\begin{equation}
    \begin{cases}
    \dot {\bar x} = A\bar x + \bar B u,
    \\
    \dot {\bar p} = -\bar p A -2p_0\alpha {\bar x}^\top.
    \end{cases}
\end{equation}
The system is subject to the initial condition $x(0)=0$, as before, but also to the terminal condition $p(T)=0$. If $p_0=0$, the terminal condition yields $p\equiv 0$. So this case is dismissed and we make the scale assumption $p_0=-1$.

Second, letting $U=[-1,1]^m$, the optimal control $u^*$ must satisfy along an optimal trajectory
\begin{equation}
    H(t,\bar x,\bar p,-1,u^*)=\max_{u\in U} H(t,\bar x,\bar p,-1,u) \quad \text{ a.e. } t\in [0,T].
\end{equation}
\textit{Notation.} $\bar{B}\in \R^{n\times m}$. We denote by $\bar B^i\in \R^n$ the $i$-th column of $\bar B$, so that, in particular, $\bar B u = \sum_{i=1}^m u_i \bar B^i$.

As a consequence, with 
$(\bar p\bar B+\psi)u=\sum_{i=1}^m(\bar p\bar B^i+\psi^i)u_i$, when $(\bar p\bar B^i+\psi^i)\neq 0$ for some $i$ in $\{1,\dots,m\}$, then the optimal control is given by $u^*_i=\mathrm{sign}(\bar p\bar B^i+\psi^i)$. These are the so-called bang arcs. On the other hand, if for some $i$ in $\{1,\dots,m\}$
\begin{equation}\label{E:sing1}
    \bar p\bar B^i+\psi^i\equiv 0
\end{equation} 
on a sub-interval of $[0,T]$, then the arc may be singular. (Arcs can be bang with respect to one control and singular with respect to another). 

Without loss of generality, consider an interval $(t_1,t_2)\subset[0,T]$ on which the first $k\leq m$ controls are singular, satisfying the condition \eqref{E:sing1}, while the remaining $m-k$ controls (if any) are of bang type (and do not switch). In this setting, we introduce constants $\varepsilon_j=\mathrm{sign}(\bar p\bar B^j+\psi^j)$ to denote the bang directions for $j=k+1,\dots ,m$.

Differentiating \eqref{E:sing1} with respect to time yields a (second) tangency condition for a trajectory to be singular:
\begin{equation}
    -\bar pA\bar B^i+2\alpha {\bar x}^\top \bar B^i + \dot{\psi}^i\equiv 0, \quad 1\leq i\leq k.
\end{equation}
Finally, differentiating once more, we have for $1\leq i\leq k$
\begin{equation}
    \bar pA^2\bar B^i
    -2\alpha {\bar x}^\top (A-A^\top)\bar B^i
    +2\alpha \sum_{j= 1}^m ({\bar B ^j})^\top{\bar B}^i u_j
    + \ddot{\psi}^i\equiv 0.
\end{equation}
Partitioning the controls into singular and non-singular leads to a linear system of $k$ equations for the singular controls $u_1^s, \dots, u_k^s$, given by, for $1\leq i\leq k$:
\begin{multline}\label{E:Gsing}
    \sum_{j= 1}^k ({\bar B ^j})^\top{\bar B}^i u_j^s
    = -\frac{1}{2\alpha} \bar pA^2\bar B^i
    + {\bar x}^\top (A-A^\top) \bar B^i
    -\frac{\ddot{\psi}^i}{2\alpha}
    \\
    -\sum_{j= k+1}^m \varepsilon_j ({\bar B ^j})^\top{\bar B}^i.
\end{multline}
Here, $\bar G\in \R^{k\times k}$ such that $\bar G_{ij}=({\bar B^j})^\top \bar B^i$ is a Gram matrix. As such, $\rank\bar G=\rank(\bar B^1,\dots ,\bar B^k)=k$, making it positive definite. From invertibility of $G$, we deduce that the singular control $(u_1^s,\dots,u_k^s)$ is fully determined by \eqref{E:Gsing} as a time-dependent affine feedback. Note that if for some $1\leq i\leq k$, some $t\in (t_1,t_2)$, equation~\ref{E:Gsing} results in $|u_i^s(\bar x,\bar  p,t)|>1$,  the singular arc cannot be sustained. Otherwise, letting $\Psi^k=(\psi^1,\dots, \psi^k)$,
\begin{multline}\label{E:feedback}
    u_s(\bar x,\bar p,t)
    =
    \bar G^{-1}\bar  B^\top\bigg(
    (A^\top-A) \bar x-\frac{1}{2\alpha}  {A^2} ^\top\bar p^\top\bigg)
    \\
    -\sum_{j= k+1}^m \varepsilon_j \bar G^{-1} \bar B^\top \bar B^j
    -
    \frac{1}{2\alpha}\bar G^{-1} \ddot{\Psi}^k(t)
\end{multline}

\subsection{Discussion of second order optimality conditions}

Let us look into second order necessary optimality conditions for the singular flow. 
See, e.g., \cite{bonnard2003singular}.
For linear control systems, Goh condition is always satisfied. 
On the other hand, Legendre-Clebsch necessary optimality condition for singular arcs in an autonomous system of  Hamiltonian $\tilde H$ is
\begin{equation}
    \frac{\partial}{\partial u}
    \frac{\partial^2}{\partial t^2}
    \frac{\partial \tilde H}{\partial u}
    \succeq 0.
\end{equation}
(See, e.g., \cite[Proposition 50]{bonnard2003singular}). 
In order to apply this criterion, we need to express our non-autonomous running cost into an autonomous one.

Consider the canonical symplectic space $T^*\R^n$, endowed with  canonical coordinates $(p,q)$. Let $H_j: T^*\R^n \to \R$, $j=0,\dots,m$, be smooth functions and let $\{\cdot,\cdot\}:T^*\R^n\to \R$ denote the Poisson bracket. We also recall the classical notation $H_{ij}=\{H_i,H_j\},H_{ijk}=\{H_i,H_{jk}\}$. 
The introduction of time as a variable is obtained considering the augmented space $T^*\R^{n+1}$, letting $\tilde p=(p,p_{n+1}), \tilde q = (q,q_{n+1})$, such that $q_{n+1}$ represents the time and $p_{n+1}$ the associated co-coordinate. For 
$1\leq j\leq m$, let
\begin{equation}
\tilde H_0(\tilde p,\tilde q) = H_0(p,q)+p_{n+1},
\tilde H_j(\tilde p,\tilde q) = H_j(p,q)+\psi^j(q_{n+1}).
\end{equation}
Then we get from direct computation.
\begin{lemma}
Letting $\tilde H_{ij}=\{\tilde H_i,\tilde H_j\},\tilde H_{ijk}=\{\tilde H_i,\tilde H_{jk}\}$, we have for $1\leq i,j\leq m$
\begin{equation}
\tilde H_{0i}(\tilde p,\tilde q)=H_{0i}(p,q)-\dot \psi^i(q_{n+1}),
\;\tilde H_{i0j}(\tilde p,\tilde q)=H_{i0j}(p,q).
\end{equation}
\end{lemma}
Since $\frac{\partial H}{\partial u_i}=p\bar B^i+\psi^i$, we can apply the above lemma with $H_i=p\bar B^i$ and $H_0=\bar pA\bar x - \alpha |\bar x|^2$. In particular, $\frac{\partial H}{\partial u_i}=\tilde H_i$.
Then for $1\leq i,j\leq m$, $\tilde H_{ij}=H_{ij}=0$, and letting $\tilde H=\tilde H_0+ \sum_{k=0}^m u_k \tilde H_k$
\begin{equation}
    \frac{\partial^2 }{\partial t^2}\frac{\partial \tilde H}{\partial u_i}
    =
    \{\tilde H,\{\tilde H , \tilde H_i\}\}
    =
    \tilde H_{00i}+\sum u_k \tilde H_{k0i}.
\end{equation}
Hence 
$
\frac{\partial }{\partial u_j}\frac{\partial^2 }{\partial t^2}\frac{\partial \tilde H}{\partial u_i}
=
H_{j0i}=2\alpha ({\bar B^j})^\top \bar B^i.
$
Consequently, along the singular trajectory where the controls $(u_i)_{1\leq i\leq k}$ are singular, the restriction of the derivative $\frac{\partial }{\partial u}\frac{\partial^2 }{\partial t^2}\frac{\partial \tilde H}{\partial u}$ to the singular components $(u_1,\dots , u_k)$ reduces to $2\alpha \bar G$. The matrix $\bar G$ being positive definite, it ensures that Legendre-Clebsch condition is satisfied.

\section{Ensemble control approach}
\subsection{Pontryagin's maximum principle statement}

For a probability measure $\nu$, we discuss a more general functional class over the space $\U=\{u\in L^\infty([0,T],\R^m) \mid \max_{1\leq i\leq m}\|u_i\|_{L^\infty}\leq 1\}$. Let $\tilde u\in \U$ to be fixed later, 
for $\eta\geq 0$ we introduce
\begin{equation}\label{E:JnuReg}
\begin{cases}
    \J_\nu^{\eta}(u)=-\langle\psi,u\rangle +\alpha \E_\nu[\langle x^\theta,x^\theta\rangle]+\frac{\eta}{2} \|u-\tilde u\|_{L^2}^2,
    \\
    \dot x^\theta =f_0^\theta (x)+\sum_{i=1}^m u_i f_i^\theta(x),
    \;
    x^\theta(0)=0,
    \;
    u\in \U.
\end{cases}
\end{equation}

Recall that $\mu_\prior$ is the prior Gaussian distribution associated with $\mathcal{N}(\theta_\prior,\Sigma_\prior)$.  In the present section, we shorten the notation to $\mu:=\mu_\prior$.  As before, denote by $x_u^\theta(\cdot)$ the solution of  $\dot x^\theta =f_0^\theta (x)+\sum_{i=1}^m u_i f_i^\theta(x)$, $x^\theta(0)=0$, $u\in \U$. We make the following assumptions.

First we assume sub-linear growth of the dynamics, in order to have existence of solutions.

\begin{assumption}\label{A:f_i}
The vector fields $(\theta,x)\mapsto f^\theta_i(x)$ are smooth over $\R^p\times \R^n$ for all $0\leq i\leq m$. Furthermore
For any $\theta\in \R^p$, there exists $c^\theta$ such that 
\begin{equation}
    \sup_{0\leq i\leq m} |f_i^\theta(x)|
    \leq 
    c^\theta (1+|x|).
\end{equation}
\end{assumption}

Similarly to the fully linear case, Assumption~\ref{A:f_i} allows to  recover Lemma~\ref{lem:conv_traj} in the nonlinear case as well, still following from,  e.g., \cite[Theorem 2.4]{scagliotti2023optimal}.
This allows to prove existence of minima for $\J_\mu$ over $\U$.

\begin{lemma}\label{T:existence_generalized}
For all $\eta\geq 0$, there exists $u\in \U$ that minimizes $\J_\mu^\eta$ in \eqref{E:JnuReg}.
\end{lemma}
\begin{pf}
First we prove that $\J_\mu$ is lower semicontinuous. 
Let $(u_N)_{N\in \N},u\in \U$, such that $u_N\to u$ in the $L^\infty$
$*$-weakly.
From Lemma~\ref{lem:conv_traj}, for any $\theta\in  \R^p$, $\sup_{t\in [0,T]}|x_{u_N}^\theta(t) -x_{u}^\theta(t)|\to 0$. We deduce that $\int_0^T |x_{u_N}^\theta(t)|^2\diff t\to \int_0^T |x_{u}^\theta(t)|^2\diff t $. From Fatou's Lemma
\begin{equation}
    \E_\mu[\langle x_{u}^\theta,x_{u}^\theta\rangle]
    \leq 
    \liminf_{N\to \infty}
    \E_\mu[ \langle x_{u_N}^\theta,x_{u_N}^\theta\rangle].
\end{equation}
The $L^2$-norm is known to be lower semi-continuous for the $L^2$-weak topology (and $L^\infty$ weak-$*$ convergence implies $L^2$-weak convergence).
Hence $\J_\mu^\varepsilon$ is lower semicontinuous. 
Since $\U$ is compact in the $L^\infty$
weak-$*$ topology, there exists a minimum for $\J_\mu^\eta$ over $\mu$.
\end{pf}

Regarding the overall behavior of the system with respect to the parameter, we assume the following.
\begin{assumption}\label{A:bound_tr}
There exists $a\in L^2(\R^p,\mu)\cap L^\infty_{loc}(\R^p,\mu)$, $b\in L^\infty_{loc}(\R^p,\mu)$ such that $|x_u^\theta(t)|\leq a(\theta)$ and $
\partial_{\theta_i}x_u^\theta(t)|\leq b(\theta)$ for all $1\leq i\leq p$, $u\in \U$, $t\in  [0,T]$, $\theta \in \R^p$.
\end{assumption}

Finally, we make an assumption on the cost, a restriction we expect to be liftable but do not address here.

\begin{assumption}\label{A:convexity}
The cost $\E_\mu \left[ \langle x^\theta_u, x^\theta_u \rangle \right]$ is  convex over $\U$. 
\end{assumption}

Under these assumptions, we recover the following extension of Pontryagin's maximum principle.

\textit{Notation.} We let $f_u^\theta=f_0^\theta+\sum_{i=1}^m u_i f_i(u)$.

\begin{theorem}\label{T:PMP}
Assume that Assumptions~\ref{A:f_i},~\ref{A:bound_tr},~\ref{A:convexity} hold. Let $H:[0,T]\times T^*\R^n\times \R \times \R^m\to \R$ be such that
\begin{multline}
    H^\theta(t,x,p,p_0,u)=
    \big\langle p, f_u^{\theta} (x)\big\rangle
    -p_0 \langle\psi(t),u\rangle +p_0\alpha \langle x,x\rangle.
\end{multline}    
If $u^*\in \U$ minimizes  $\J_\mu:=\J_\mu^0$ in \eqref{E:JnuReg}
there exists $(x,p):\R^p\to C^1([0,T],T^*\R^{n})$ and a scalar $p_0\leq 0$, with $(p,p_0)\neq 0$, such that for a.e. $t\in [0,T]$, for all $\theta\in \R^p$, 
\begin{equation}
\begin{aligned}
\dot{x}^\theta(t) =& \frac{\partial H}{\partial p}(t,x^\theta(t),p^\theta(t),p_0,u^*(t)),
\\
\dot{p}^\theta(t) = &-\frac{\partial H}{\partial x}(t,x^\theta(t),p^\theta(t),p_0,u^*(t)),
\end{aligned}
\end{equation}
with boundary conditions
$x(0)= 0$,  $p(T) = 0$.
Furthermore, for almost every $t\in [0,T]$,
\begin{multline}
\max_{u\in [-1,1]^m}\E_{\mu}\left[H^\theta(t,x^\theta(t),p^\theta(t),p_0,u)\right]
\leq \\
\E_\mu\left[H^\theta(t,x^\theta(t),p^\theta(t),p_0,u^*(t))\right].
\end{multline}
\end{theorem}

Before discussing the proof, let us just remark that our case of interest is well covered by this statement.

\begin{proposition}
Systems of the form $\dot x^\theta = A^\theta x + B^\theta u$, with $A^\theta=A_0+\sum \theta_i A_i$ and $B^\theta=B_0+\sum \theta_i B_i$ satisfy Assumptions~\ref{A:f_i},~\ref{A:bound_tr},~\ref{A:convexity}.
\end{proposition}
\begin{pf}
Assumptions~\ref{A:f_i} is immediate. Assumptions~\ref{A:bound_tr} holds because 
$\mu$ is a Gaussian distribution over $\R^p$. Indeed there exist constants $\alpha,\beta>0$ such that, for every $u \in \mathcal{U}$,
\begin{equation}
    |x_u^\theta(t)|\leq  \e^{\alpha |\theta| }\beta |\theta|=:a(\theta) \in L^2(\R^p,\mu)\cap L^\infty_{loc}(\R^p,\mu).
\end{equation}
On the other hand, for $1\leq i\leq p$, if in general $x^\theta$ is such that $\dot x^\theta =f_0^\theta (x^\theta)+\sum_{j=1}^m u_j f_j^\theta(x^\theta)=f_u^\theta(x^\theta)$, then (the vector fields being smooth w.r.t. $\theta$)
$\partial_{\theta_i} \dot{x}^\theta
=
\partial_{\theta_i}f_u^\theta (x^\theta)+ \nabla_x f_u^\theta (x^\theta) \partial_{\theta_i}x^\theta$
and $\partial_{\theta_i} \dot{x}^\theta(0)=0$. Hence $\max_{1\leq i\leq p}|\partial_{\theta_i}x_u^\theta(t)|$ is a well defined quantity.
Letting
\begin{equation}
\tilde{A}^\theta=
\begin{pmatrix}
    A^\theta&&&(0)\\
    A_1 &\ddots\\
    \vdots&(0)&\;\ddots \\
    A_p&& &A^\theta
\end{pmatrix},
\;
\tilde B^\theta =\begin{pmatrix}
    B^\theta\\
    B_1\\
    \vdots \\
    B_p
\end{pmatrix},
\;
\tilde x^\theta
=
\begin{pmatrix}
    x^\theta\\
    \partial_{\theta_1}x^\theta\\
    \vdots \\
    \partial_{\theta_p}x^\theta
\end{pmatrix},
\end{equation}
we have again 
$\dot{\tilde x} = \tilde A^\theta \tilde x + \tilde B^\theta u$. 
In particular, there exist constants $\tilde \alpha,\tilde \beta>0$ such that, for every $u \in \mathcal{U}$,
\begin{equation}
    |\tilde x_u^\theta(t)|\leq  \e^{\tilde \alpha |\theta| }\tilde \beta |\theta|=:b(\theta) \in L^\infty_{loc}(\R^p,\mu).
\end{equation}

Assumption~\ref{A:convexity} follows from the linearity in control of
$x_{u}^\theta(t)=\int_0^t \e^{A^\theta(t-s)}B^\theta u(s)\diff s$. 
For any two $u,v\in \U$, the map $\phi_\theta:\lambda\mapsto \langle x_{u+\lambda v}^\theta, x_{u+\lambda v}^\theta\rangle_{L^2}$ satisfies
$
\phi_\theta''(\lambda)
=
2 \langle x_{v}^\theta, x_{v}^\theta\rangle_{L^2}.
$
That is, as soon as $v\neq 0$, $\phi_\theta''(\lambda)>0$, proving the (strict) convexity of 
$\U \ni u \mapsto \E_\mu \left[ \langle x^\theta_u, x^\theta_u \rangle \right]$.
\end{pf}

The proof of Theorem~\ref{T:PMP} is achieved by showing that the ensemble problem is the $\Gamma$-limit of a collection of finite dimensional optimal control problems.
Recall (see, e.g., \cite[Proposition 8.1]{dal2012introduction}) that a sequence $(F_N)_{N\in \N}$ of functionals over a first-countable topological space $E$  $\Gamma$-converges to $F$ if and only if the following two conditions hold
\begin{enumerate}[label=(\roman*)]
    \item\label{C:Gamma1} For any $f\in E$, and $(f_N)_{N\in \N}\in E$ such that $f_N\to f$,
    \begin{equation}
        F(f)\leq \liminf_{N\to \infty} F_N(f_N).
    \end{equation}
    \item\label{C:Gamma2} For every $f\in E$, there exists a sequence $(f_N)_{N\in \N}\in E$,  $f_N\to f$, such that 
        $F_N(f_N)\to F(f)$.
\end{enumerate}


Similarly to \cite{scagliotti2023optimal}, we approximate $\mu$ by a measure with finite support. Such measures are dense in the set of measures with finite first moment in the $W^1$ topology\footnote{
Weak convergence 
against
1-Lipschitz-continuous test functions.}. See, e.g., \cite[Theorem 2.2.7]{Panaretos2020}. Assumption~\ref{A:bound_tr} implies the following.

\begin{lemma}\label{L:approxMu}
There exists a sequence $(\mu_N)_{N\in \N}$ such that 
\begin{enumerate}
    \item For any integer $N>0$, $\mu_N$ is a probability measure on $\R^p$ supported on a finite collection of points within $\bar B(0,N)$, the closed ball of radius $N$ centered at $0$.
    \item The sequence $\mu_N$ converges to $\mu$ in the $W^1$ topology.
    \item For any $\eta\geq 0$ and any sequence $u_N\in \U$,
    $\J_{\mu_N}^\eta(u_N)-\J_\mu^\eta(u_N)\to 0$.
\end{enumerate}

\end{lemma}

\begin{proposition} \label{prop:G-conv}
{  For all $\eta\geq 0$,} the sequence $(\J_{\mu_N}^\eta)_{N\in \N}$ $\Gamma$-converges towards $\J_\mu^\eta$ in the space $\U$ endowed with the $L^\infty$ weak-$*$ topology.
\end{proposition}

\begin{pf}
We start with point~\ref{C:Gamma1}. Let $(u_N)_{N\in \N}$ in $\U$, $u\in L^\infty([0,T],\R^m)$ such that that $u_N\to u$ for the  $L^\infty$ weak-$*$ convergence. 
Since $\U$ is convex and strongly closed, it is weakly closed as well, so that $u$ belongs to $\U$.
From Lemma~\ref{L:approxMu}, we have $\J_{\mu_N}^\eta(u_N)-\J_\mu^\eta(u_N)\to 0$.
Second, from Lemma~\ref{lem:conv_traj}, we have the pointwise convergence of $|x_{u_N}^\theta(t)|^2$ towards $|x_{u}^\theta(t)|^2$ for all $(t,\theta)\in [0,T]\times \R^p$.  Since $|x_{u_N}^\theta(t)|^2\leq a(\theta)^2$ and $a^2\in L^1([0,T]\times \R^p, \diff t \otimes \mu)$, the dominated convergence theorem applies to prove 
$\J_\mu^\eta(u_N)\to \J_\mu^\eta(u)$.
As a result $\J_{\mu_N}^\eta (u_N) - \J_\mu^\eta(u)\to 0$ and point~\ref{C:Gamma1} holds.

Point~\ref{C:Gamma2} follows from Lemma~\ref{L:approxMu}: pick $u\in\U$, and the sequence $u_N=u$, $N\in \N$; then  $\J_\mu^\eta(u)=\lim_{N\to \infty} \J_{\mu_N}^\eta(u)$.
\end{pf}

Applying standard $\Gamma$-convergence arguments (based on lower semi-continuity of each element of $\J_\mu^\eta$), we get the following.
\begin{corollary} \label{C:conclusion}
Let $\eta\geq 0$.
Any sequence $(u^*_N)_{N\in \N}$ in $\U$ such that $u^*_N \in \arg \min_{\mathcal U} \J_{\mu_N}^\eta$ is pre-compact in the $L^ \infty$ weak-$*$ topology, and any accumulation point is a minimizer for $\J_\mu^\eta$.
Moreover, we have that
\begin{equation} \label{eq:conv_minima}
    \lim_{N\to \infty} \min_\U \J_{\mu_N}^\eta = \min_\U \J_\mu^\eta,
\end{equation}
and, if $u^*_N \to u^*$ in the $L^\infty$ weak-$*$ topology, then
\begin{equation}
\lim_{N\to \infty} \E_{\mu_N}  \left[ \langle x_{u^*_N}^\theta , x_{u^*_N}^\theta \rangle \right] 
= \E_{\mu} \left[\langle x_{u^*}^\theta , x_{u^*}^\theta \rangle \right],
\end{equation}
\begin{equation}\label{E:conv_reg}
    \lim_{N\to \infty} \frac{\eta}{2}\| u^*_N- \tilde u\|^2_{L^2} = \frac{\eta}{2}\| u^*- \tilde u\|^2_{L^2}.
\end{equation}
\end{corollary}

Finally, we have the following classical property thanks to the Radon--Riesz property of $L^2(0,T)$.

\begin{corollary} \label{C:conclusion_ref} 
Let $\eta>0$. Under Assumption~\ref{A:convexity}, any sequence 
$(u^*_N)_{N\in \N}$ such that $u^*_N \in \arg \min_{\mathcal U} \J_{\mu_N}^\eta$ is convergent in the \emph{$L^2$-norm} to the unique minimizer of $\J_\mu^\eta$. 
\end{corollary}


We are now  ready for the proof of the PMP extension.

\begin{pf}(Proof of Theorem~\ref{T:PMP}.) 
From Lemma~\ref{T:existence_generalized}, we know that at least one minimizer, that we denote $\uopt$, exists. We proceed
to apply the above results using $\tilde u = \uopt$. It is immediate to see that for all $\eta>0$, for all $u\in \U$
\begin{equation}
\J_{\mu}^\eta(u)\geq \J_{\mu}^0(u) \geq \inf_{v\in \U}\J_{\mu}^0(v)= \J_{\mu}^0(\uopt) = \J_{\mu}^\eta(\uopt).
\end{equation}
As a result, the minimizer of $\J_{\mu}^\eta$ in Corollary~\ref{C:conclusion_ref} is $\uopt$.

The next step is to apply Pontryagin's maximum principle for the cost function $\J_{\mu_N}^\eta$ with a fixed value of $\eta>0$ and $N\in \N$. Since $\mu_N$ is atomic, let $(\theta^k_N)_{1\leq k\leq k_N}$ denote the (finite) collection of parameters such that $\mu_N(\theta^k_N)\neq 0$. Then we can describe the state in the optimal control problem as the collection $x:=(x^{\theta_1},\dots,x^{\theta_{k_N}}) \in (\R^{n})^{k_N}$. For $p:=(p^{\theta_1},\dots,p^{\theta_{k_N}}) \in T^*(\R^{n})^{k_N}\simeq (T^*\R^{n})^{k_N}$,
the associated Hamiltonian is, for almost every $t\in[0,T]$,
\begin{multline}
    H_N^\eta(t,x,p,p_0,u)=
    -p_0\langle \psi(t), u\rangle
    +p_0\frac{\eta}{2}| u-\uopt(t)|^2+
    \\
    \sum_{k=1}^{k_N}
    \left(\big\langle p^{\theta^k_N}, f_u^{\theta^k_N} (x^{\theta^k_N})\big\rangle
    +
    p_0\alpha \mu_N(\theta_N^k)\langle x^{\theta^k_N},x^{\theta^k_N}\rangle\right).
\end{multline}
In particular, we can write $(p,p_0)\neq 0$, the evolution
\begin{equation}
\frac{\diff p^{\theta_N^k}}{\diff t} 
=  
-\frac{\partial}{\partial x^{\theta^k_N}}
\big\langle p^{\theta^k_N}, f_u^{\theta^k_N} (x^{\theta^k_N})\big\rangle
+2p_0\alpha
\mu_N(\theta_N^k)x^{\theta^k_N},
\end{equation}
the endpoint transversality condition $p(T)=0$, and that an optimal control $u_N^*$ satisfies along its associated trajectory $(x,p)$, for $\text{a.e. } t\in [0,T]$, and all $u\in [-1,1]^m$:
\begin{equation}\label{E:PMP}
H_N^\eta(t,x(t),p(t),p_0,u)\leq H_N^\eta(t,x(t),p(t),p_0,u_N^*(t)).
\end{equation}

We now change variables.
For all $\theta\in \R^p$, define $q^{\theta}$ by
\begin{equation}\label{E:dynP}
\frac{\diff q^{\theta}}{\diff t} 
=  
-\frac{\partial}{\partial x^{\theta}}
\big\langle q^{\theta}, f_u^{\theta} (x^{\theta})\big\rangle
+2p_0\alpha
x^{\theta}, \quad  q^{\theta}(T)=0.
\end{equation}
Take
$p^{\theta_N^k}=\mu_N(\theta_N^k)q^{\theta_N^k}$
and write $H_N^\eta$ as a map of $(x,q)$:
\begin{multline}
    H_N^\eta(t,x,q,p_0,u)=
    -p_0\langle \psi(t), u\rangle
    +p_0\frac{\eta}{2}| u-\uopt(t)|^2+
    \\
    \sum_{k=1}^{k_N}
    \mu_N(\theta_N^k)
    \left(\big\langle q^{\theta^k_N}, f_u^{\theta^k_N} (x^{\theta^k_N})\big\rangle
    +
    p_0\alpha \langle x^{\theta^k_N},x^{\theta^k_N}\rangle\right).
\end{multline}
We can recognize
\begin{multline}
H_N^\eta(t,x,q,p_0,u)
=
\E_{\mu_N}\left[H^\theta(t,x^\theta(t),q^\theta(t),p_0,u)\right] 
\\+ p_0 \frac{\eta}{2}| u-\uopt(t)|^2.
\end{multline}

Now consider a sequence of optimal controls $(u^*_N)_{N\in \N}$ in $\U$ such that $u^*_N \in \arg \min_{\mathcal U} \J_{\mu_N}$.
Owing to Corollary ~\ref{C:conclusion_ref} we have that $\| u^*_N - \uopt\|_{L^2} \to 0$. Hence, we can extract a (non-relabeled) subsequence such that $u^*_N(t) \to u^*(t)$ for a.e.~$t\in [0,T]$. 
From Lemma~\ref{lem:conv_traj}, we know that as $u^*_N\to u^*$, the trajectory $x_{u_N^*}^{\theta}$ converges uniformly in time to $x_{u^*}^{\theta}$ for every $\theta$.
A similar uniform convergence statement can be obtained for the joint  state and covector trajectory. 
In particular, as in Corollary~\ref{C:conclusion}, we have $\E_{\mu_N} \left[ \langle x_{u^*_N}^\theta , x_{u^*_N}^\theta \rangle \right] \to \E_{\mu} \left[ \langle x_{u^*}^\theta , x_{u^*}^\theta \rangle \right]$.
Finally, $\langle\psi (t),u_N^*(t) \rangle \to \langle \psi(t) ,u^*(t) \rangle $ for a.e.~$t\in [0,T]$. 
Hence, we have for a.e.~$t\in [0,T]$
\begin{multline}
    \lim_{N\to \infty} 
    H_N^\eta(t,x(t),q(t),p_0,u^*_N(t)) =\\
     \E_\mu\left[H^\theta(t,x^\theta(t),p^\theta(t),p_0,\uopt(t)\right].
\end{multline}
Then considering the set of times for which \eqref{E:PMP} holds, and intersecting over all $N$, we can also pass to the limit in the inequality to state, for a.e. $t\in [0,T]$, for all  $u\in [-1,1]^m$, 
\begin{multline}\label{E:PMP2}
\E_\mu\left[H^\theta(t,x^\theta(t),p^\theta(t),p_0,u\right]
 + p_0 \frac{\eta}{2}| u-\uopt(t)|^2
\leq 
\\
\E_\mu\left[H^\theta(t,x^\theta(t),p^\theta(t),p_0,\uopt(t)\right].
\end{multline}
(Recall $p_0\leq 0$.) This equation has been obtained for an arbitrary value of $\eta$. Hence considering the set of times for which \eqref{E:PMP2} holds with the values $\eta_k=1/k$, and intersecting of all $k\geq 1$, we can pass to the limit to get for a.e. $t\in [0,T]$, for all  $u\in [-1,1]^m$,
\begin{multline}\label{E:PMP3}
\E_\mu\left[H^\theta(t,x^\theta(t),p^\theta(t),p_0,u\right]
\leq
\\
\E_\mu\left[H^\theta(t,x^\theta(t),p^\theta(t),p_0,\uopt(t)\right].
\end{multline}
This finishes the proof of the statement.
\end{pf}

\subsection{Ensemble case}
Now we propose to apply the generalized Pontryagin's maximum principle to our case. We consider a distributed state $(x^\theta)_{\theta\in \R^p}$ and covector $(p^\theta)_{\theta\in \R^p}$, such that $(p^\theta,x^\theta)\in T^*\R^n$ for all $\theta$. The pointwise Hamiltonian takes the form
\begin{equation}\label{E:Ham1bis}
    H^\theta(t,x^\theta,p^\theta,p_0,u) = p^\theta Ax^\theta  + p_0 \alpha |x^\theta|^2 + (p^\theta B^\theta-p_0\psi)u.
\end{equation}
Here, $p_0\leq 0$ is constant such that $\left((p^\theta)_{\theta\in \R^p},p_0\right)$ is nontrivial. 
The dynamics of the distributed state-covector system follow from the  Hamiltonian trajectory
\begin{equation}
    \begin{cases}
    \dot x^\theta = Ax^\theta + B^\theta u,
    \\
    \dot p^\theta = -p^\theta A -2p_0\alpha (x^\theta)^\top.
    \end{cases}
\end{equation}
The system is subject to the boundary conditions $x^\theta(0)=0$, $p^\theta(T)=0$. If $p_0=0$ the terminal condition yields that $p^\theta\equiv 0$ for all $\theta\in  \R^p$ as well. So this case can be dismissed and we make the scale assumption $p_0=-1$. The rest of the proof follows a similar path to Section~\ref{S:classical}.

With $U=[-1,1]^m$, the optimal control $u^*$ must satisfy along an optimal trajectory, for all $u\in U$,  a.e. $t\in [0,T]$,
\begin{equation}
    \E_{\mu}\left[H^\theta(t,x^\theta,p^\theta,-1,u^*(t))\right]\geq \E_{\mu}\left[H^\theta(t,x^\theta,p^\theta,-1,u)\right].
\end{equation}
Since 
$(p^\theta B^\theta+\psi)u=\sum_{i=1}^m(p^\theta (B^\theta)^i+\psi^i)u_i$, when $(\E_\mu\left[p^\theta (B^\theta)^i\right]+\psi^i)\neq 0$ for some $i$ in $\{1,\dots,m\}$, then the optimal control is given by $u^*_i=\mathrm{sign}(\E_\mu\left[p^\theta (B^\theta)^i\right]+\psi^i)$. 
We again find bang arcs.
If for some $i$ in $\{1,\dots,m\}$
\begin{equation}\label{E:sing1bis}
    \E_\mu\left[p^\theta (B^\theta)^i\right]+\psi^i\equiv 0
\end{equation} 
on a sub-interval of $[0,T]$, then the arc may be singular.

Again, consider $(t_1,t_2)\subset[0,T]$ on which the first $k\leq m$ controls are singular, satisfying \eqref{E:sing1bis}, while the remaining $m-k$ controls (if any) are of bang type (and do not switch). Let $\varepsilon_j=\mathrm{sign}(\E_\mu\left[p^\theta (B^\theta)^j\right]+\psi^j)\in \{\pm 1\}$ denote the bang directions for $j\in\{k+1,\dots ,m\}$.

Differentiating \eqref{E:sing1bis} yields the tangency condition:
\begin{equation}
    -\E_\mu\left[p^\theta A(B^\theta)^i\right]+2\alpha \E_\mu\left[(x^\theta)^\top (B^\theta)^i\right] + \dot{\psi}^i\equiv 0, \; 1\leq i\leq k.
\end{equation}
Differentiating once more,  for all $1\leq i\leq k$, 
\begin{multline}
    \E_\mu\left[p^\theta A^2(B^\theta)^i\right]
    +2\alpha \E_\mu\left[(x^\theta)^\top A^\top (B^\theta)^i-(x^\theta)^\top A(B^\theta)^i\right]
    \\
    +2\alpha \sum_{j= 1}^m \E_\mu\left[{(B^\theta)^j}^\top (B^\theta)^i\right] u_j
    + \ddot{\psi}^i\equiv 0.
\end{multline}
We follow the same strategy as Section~\ref{S:classical}, partitioning between singular controls $u_1^s, \dots, u_k^s$ and bang arcs. We need to account for the expectation in order to find a feedback form.
Let us recall a notation introduced early in the paper. $(B^\theta)^j=B_0^j+B^j_1\theta_1+\dots +B^j_p\theta_p$, so we set  $B^j\in \R^{n\times p}$ to be such that $B^j\theta = B^j_1\theta_1+\dots +B^j_p\theta_p$. Then,
\begin{equation}
\begin{aligned}
\E_\mu\left[{(B^\theta)^j}^\top (B^\theta)^i\right]
&=
\E_\mu\left[\theta^\top  {B^j}^\top  B^i \theta\right]
+
\theta_\prior^\top  {B^j}^\top  B^i_0
\\
&\quad 
+
{B^j_0}^\top  B^i \theta_\prior
+
{B^j_0}^\top  B^i_0.
\end{aligned}
\end{equation}
Usual manipulations allow to show
\begin{equation}
\E_\mu\left[\theta^\top  {B^j}^\top  B^i \theta\right]
=
\Tr \left( {B^j}^\top B^i \Sigma_\prior\right)
+
\theta_\prior^\top {B^j}^\top B^i  \theta_\prior,
\end{equation}
in order to recover
\begin{equation}
\E_\mu\left[{(B^\theta)^j}^\top (B^\theta)^i\right]
=
\Tr \left( {B^j}^\top B^i \Sigma_\prior\right)
+
(\bar B^j)^\top \bar B^i.
\end{equation}
One can recognize in $(\bar B^j)^\top \bar B^i={(B^{\theta_\prior})^j}^\top
(B^{\theta_\prior})^i$, the $(i,j)$-coefficient of the Gram matrix $\bar G$ introduced in Section~\ref{S:classical}. 
Regarding the term $\Tr \left( {B^j}^\top  B^i \Sigma_\prior\right)$, we can check that $G^{\Sigma_\prior}\in \R^{k\times k}$ given by $G^{\Sigma_\prior}_{ij}=\Tr \left( {B^j}^\top  B^i \Sigma_\prior\right)$ is again a Gram matrix, this time for the family of elements $\left(B^i\Sigma_\prior^{1/2}\right)_{1\leq i\leq m}$, with respect to the Frobenius inner product on $n\times p$ matrices. As a result, assuming  $\Sigma_\prior$ is not degenerate, $\rank G^{\Sigma_\prior}=\rank (  B^1,\dots ,  B^k)$.
Under Assumption~\ref{A:Btheta}, the family $\left(  B^i\right)_{1\leq i\leq m}$ is free, so $G^{\Sigma_\prior}$ is positive definite.
As in the classical case, 
we deduce that the singular control $(u_1^s,\dots,u_k^s)$ is fully determined as a time-dependent affine feedback. Again, the singular arc cannot be sustained if it results in $|u_i^s(x,p,t)|>1$, for some $1\leq i\leq k$, $t\in (t_1,t_2)$. Otherwise, as in \eqref{E:feedback}, letting $\Psi^k=(\psi^1,\dots, \psi^k)$
\begin{multline}
    u_s(x,p,t)
    =
    G^{-1}
    \E_\mu\!\left[ {B^\theta}^\top\!
    \left(
    (A^\top-A) x^\theta 
    -
    {A^2}^\top 
    {p^\theta}^\top /2\alpha\right)\right]
    \\
    -\sum_{j= k+1}^m \varepsilon_j G^{-1}\E_\mu\left[{(B^\theta)}^\top (B^\theta)^j\right]
    -G^{-1}
    \frac{\ddot{\Psi}^k(t)}{2\alpha}
\end{multline}

It is not yet known how to extend second order arguments to ensemble control, but the very similar nature of the optimal flow leads naturally to such questions, which are left to future endeavors.

\section{Numerical experiments}
\begin{figure}[hb]
    \centering
\includegraphics[width=\linewidth]{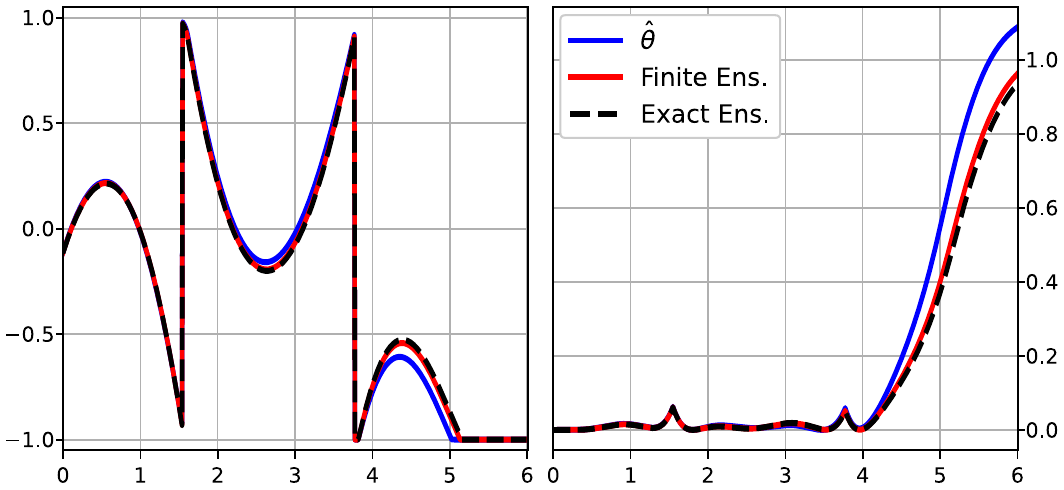}
    \caption{Left: the computed controls $u_{\hat \theta}$ (blue), with $u_{\mu_{N}}$ (red) and $u_{\mu_{\mathrm{prior}}}$ (black). Right: the graph of $t\mapsto |x^{\theta_{\mathrm{true}}}_u(t)|^2$ 
    under each of these policies.}
    \label{figure_1}
\end{figure}

We consider a damped linear harmonic oscillator where we induce an acceleration tuned by the control $u\in [-1,1]$:
\begin{equation*}
    \left(
    \begin{matrix}
        \dot x\\
        \dot v
    \end{matrix}
    \right) =
        \left(
    \begin{matrix}
         0 & 1\\
         -\kappa & -\delta
    \end{matrix}
    \right)
        \left(
    \begin{matrix}
         x\\
         v
    \end{matrix}
    \right)
    +       \left[ \left(
    \begin{matrix}
         0\\
         1
    \end{matrix}
    \right)
    + \left(
    \begin{matrix}
         0\\
         \theta
    \end{matrix}
    \right)
    \right] u,
    \quad C = (\begin{matrix}
        0 & 1
    \end{matrix}),
\end{equation*}
over time $[0,T]$, with $\kappa = 2, \delta =0.25$, and where $\theta$ is the unknown parameter with a prior distribution $\mu_{\mathrm{prior}}\sim \mathcal{N}(0,0.5)$. 
We aim to design an experiment to update $\mu_{\mathrm{prior}}$, and we want the produced trajectory $t\mapsto x_u(t)$  to satisfy $|x_u(t)|^2\leq 1$ for every $t\in [0,T]$.
After setting $\alpha=1.2$, we solved numerically the optimal control problems stated in \eqref{eq:problemClassic}, corresponding to the single-system case $\bar \theta= \hat \theta= 0$, and \eqref{eq:problemEnsemble}, 
yielding controls $u_{\hat \theta}$ and $u_{\mu_{\mathrm{prior}}}$, respectively.
Following Remark~\ref{rmk:exact_ensemble_comput}, we reformulated \eqref{eq:problemEnsemble} as a single-system problem.
From the $\Gamma$-convergence established in Proposition~\ref{prop:G-conv}, we also solved \eqref{eq:problemEnsemble} for $\mu_N$ an approximation of $\mu_\prior$ consisting of $N=51$ equi-spaced atoms, giving $u_{\mu_N}$.
We used the automatic differentiation tools of PyTorch, and we performed $10^3$ steps of the gradient descent with step-size $h=0.5$, using as initial guess the explicit solution of the optimal control problem corresponding to $\alpha=0$.
The results are shown in Fig.~\ref{figure_1}.
We observe that the computed controls are qualitatively close to each other. However, when we run the experiment designed by the computed controls on the `true' system corresponding to $\theta = \theta_{\mathrm{true}}=0.25$, we observe that the trajectory driven by $u_{\hat \theta}$ violates the constraint on the trajectory magnitude.
To complete the Bayesian posterior updates, for $u_{\hat \theta}$, $u_{\mu_{N}}$ and $u_{\mu_{\mathrm{prior}}}$ we constructed the measurement $Y_{\mathrm{avg}}$ as prescribed in \eqref{eq:Y_signal}, with $\sigma = 0.25$. 
Following Section~\ref{S:Bayesian}, we computed the posterior distributions corresponding to the experiment result, getting:
\begin{itemize}
    \item For $u_{\hat \theta}$, $\mu_{\mathrm{post}}=0.1138$ and $\Sigma_{\mathrm{post}}=0.2775$;
    \item For $u_{\mu_{N}}$, $\mu_{\mathrm{post}}=0.1034$ and $\Sigma_{\mathrm{post}}=0.2981$;
    \item For $u_{\mu_{prior}}$, $\mu_{\mathrm{post}}=0.1009$ and $\Sigma_{\mathrm{post}}=0.3031$.
\end{itemize}
The results above should be compared with $\mu_{\mathrm{prior}}=0$ and $\Sigma_{\mathrm{prior}}=0.5$, and with $\theta_{\mathrm{true}}=0.25$.

\section{Conclusion}

We introduced an optimal-control formulation for Bayesian input design and compared two viewpoints: a classical approach based on a single nominal parameter and an ensemble formulation that optimizes the prior-averaged cost. The ensemble viewpoint naturally incorporates prior uncertainty and led us to prove a tailored Pontryagin maximum principle for the Gaussian linear setting, including the unbounded-parameter regime treated in this paper, and to characterize the resulting bang–bang/singular flow. Future work will pursue several complementary directions: further generalize the ensemble PMP to broader functional settings, remove the linear dependence on the parameter and deepen the conceptual links between optimal ensemble control and experiment design.

\bibliography{biblio}

\end{document}